\documentclass[11pt]{article}

\usepackage[pdftex]{graphicx}
\usepackage{amssymb}
\usepackage{amsthm}

\makeatletter
\renewcommand\section{\@startsection{section}{1}{\z@}%
                                  {-3.5ex \@plus -1ex \@minus -.2ex}%
                                  {2.3ex \@plus.2ex}%
                                  {\normalfont\large\bfseries}}
\makeatother

\DeclareGraphicsExtensions{.jpg}

\begin{document}

%\begin{frontmatter}



\title{Transposition diameter on circular binary strings}


\author{Misa Nakanishi \thanks{E-mail address : nakanishi@2004.jukuin.keio.ac.jp}}
\date{}
\maketitle

%\footnote{E-mail address : nakanishi@2004.jukuin.keio.ac.jp}

%\address{}

\begin{abstract}
On the string of finite length, a (genomic) transposition is defined as the operation of exchanging two consecutive substrings. The minimum number of transpositions needed to transform one into the other is the transposition distance, that has been researched in recent years. In this paper, we study transposition distances on circular binary strings. A circular binary string is the string that consists of symbols $0$ and $1$ and regards its circular shifts as equivalent. The property of transpositions which partition strings is observed. A lower bound on the transposition distance is represented in terms of partitions. An upper bound on the transposition distance follows covering of the set of partitions. The transposition diameter is given with a necessary and sufficient condition. \\
keyword : strings, sorting by transpositions, partitions, NP-hard problems
\end{abstract}

%\end{frontmatter}

%\linenumbers

\section{Introduction}
\ \ \ \ \ Let $A = \{0, 1, \dots, a-1\}$ be a finite set, called an alphabet. For the string of finite length, which is $S = S(1) S(2) \dots S(n)$ $(S(i) \in A \ (1 \leq i \leq n))$, a transposition $\mathcal{T}_{ijk}$ is defined as the operation of exchanging two consecutive substrings $S(i) S(i+1) \dots S(j-1)$ and $S(j) S(j+1) \dots S(k-1)$ $(1 \leq i < j < k \leq n+1)$.

For two strings $S$ and $T$ that have the same multiplicity of symbols, the minimum number of transpositions needed to transform $S$ into $T$ is called the transposition distance.

A transposition is studied as one of operations on genome rearrangement~\cite{B. Alberts}. On the other hand, a transposition is concerned with a combinatorial problem such as sorting a bridge hand~\cite{H. Eriksson}. Both of the studies have the object to find the minimum transposition sequence which transforms one into the other.

In the 1990s, transpositions on permutations were investigated. Bafna and Pevzner~\cite{V. Bafna} gave a lower bound and an upper bound on the transposition distance between given pairs of permutations, and the transposition diameter, which is the maximum transposition distance taken over all pairs of permutations. Afterwards, Eriksson et al.~\cite{H. Eriksson} improved the upper bound and introduced circular permutations they call toric permutations. Hartman~\cite{T. Hartman} showed circular permutations correspond to linear ones of which the length reduces by 1.

Then, transpositions on strings, which can contain redundant symbols, were researched as an expansion. Christie and Irving~\cite{D. A. Christie} determined the transposition diameter and its string pairs on binary strings. Radcliffe, Scott and Wilmer~\cite{A. J. Radcliffe} led NP-hardness of transposition distances on binary strings from 3-PARTITION~\cite{Computers and Intractibility}. Therefore, on arbitrary strings, the complexity is also NP-hard.

In this paper, we study transpositions on circular binary strings. A circular binary string is the string that consists of symbols $0$ and $1$ and regards its circular shifts as equivalent. In Section 2, terminology and definitions are introduced. In Section 3, the property of transpositions is observed, which partition
strings. A partition is constituted by parts between $1$'s where weights as the number of $0$'s belong.
Several functions are introduced that sum up the minimal or maximal weights and give a lower bound on
transposition distances. A consecutive part is covered by finite number of
parts. In this regard, an upper bound on transposition distances follows. Finally, 
the transposition diameter is represented by a necessary and sufficient condition. The distance between two strings is determined in polynomial time when it is a diameter.
%\label{}

\section{Preliminary}
\ \ \ \ \ For an alphabet $A = \{0, 1\}$, a binary string $S$ of length $n$, that is $S \in A^{n}$, is denoted by
\[S = S(1) S(2) \dots S(n) \ in\ which\ S(i) \in A \ (1 \leq i \leq n).\]

For a binary string $S$, a \textit{transposition} $\mathcal{T}_{ijk} \ (1 \leq i < j < k \leq n+1)$ is marked down as
\[S(1) \dots S(i-1)[S(i) \dots S(j-1)][S(j) \dots S(k-1)]S(k) \dots S(n),\]
that is transformed into
\[S(1) \dots S(i-1) S(j) \dots S(k-1) S(i) \dots S(j-1) S(k) \dots S(n).\]

In the case that binary strings $S$ and $T$ have the same multiplicity of symbols, the \textit{transposition distance} is defined as the minimum number of transpositions needed to transform $S$ into $T$, that is denoted by $d(S, T)$. In this paper, transposition distances are defined to such string pairs.

Let $0^{i}$ be consecutive run of $0$ of length $i$ and let $0^{\ast}$ be consecutive run of $0$ of length 0 or more. For a string $S$, $S^{i}$ (or $(S)^{i}$) represents the string that concatenates $S$ $i$ times. $S \cdot T$ or simply $ST$ represents the concatenation of $S$ and $T$. As the symbol of a concatenation of strings, we use $\Sigma$. For example, $\Sigma_{i = 1}^{3}(0^{i}1) = 010010001$. In this paper, it is assumed that the number of $1$'s is less than or equal to the number of $0$'s on a binary string, otherwise the symbols are interchanged.

A \textit{circular} binary string sees circular shifts of the string as equivalent. A equivalence relation $\sim$ is defined as
\[S(1) \dots S(n) \sim S(i) S(i+1) \dots S(n) S(1) \dots S(i-1) \ (1 \leq i \leq n).\]
A circular binary string $S$ as an equivalence type is denoted by
\[S = (S(1) S(2) \dots S(n)).\]

A transposition for a circular binary string is defined in the same manner
as a transposition for a linear string. Note that it separates a circular binary string into 3 segments and connects them in different order.

In this paper, all variable numbers are within range of the natural numbers.

\section{Transposition distances on circular binary strings}
\ \ \ \ \ A pair of circular binary strings is represented by
\[S = (0^{s_{1}}1 \dots 0^{s_{k}}1) \ and\ T = (0^{t_{1}}1 \dots 0^{t_{k}}1).\]

A transposition is classified into type (T1) or type (T2) in terms of the number of parts it operates on. The index of each part is taken as modulo $k$.

\begin{enumerate}

\item[(T1)]
The transposition of this type operates on 3 parts. \\
For parts $\alpha, \beta, \gamma \in \{1, \dots, k\}$ $(\alpha < \beta < \gamma)$ on $S$
, and $0 \leq x \leq s_{\alpha}$, $0 \leq y \leq s_{\beta}$ and $0 \leq z \leq s_{\gamma}$,
\[ S = (10^{s_{\alpha} - x}[0^{x}(10^{\ast})^{\beta - \alpha - 1}10^{s_{\beta} - y}][0^{y}(10^{\ast})^{\gamma - \beta - 1}10^{s_{\gamma} - z}]0^{z}(10^{\ast})^{\alpha - \gamma - 1 + k}) \]
is transformed into
\[ S' = (10^{s_{\alpha} - x + y}(10^{\ast})^{\gamma - \beta - 1}10^{s_{\gamma} - z + x}(10^{\ast})^{\beta - \alpha - 1}10^{s_{\beta} - y + z}(10^{\ast})^{\alpha - \gamma - 1 + k}). \]
If there exists $\alpha', \beta', \gamma' \in \{1, \dots, k\}$ on $T$ such that 
\[s_{\alpha} - x + y = t_{\alpha'},\ s_{\gamma} - z + x = t_{\beta'} \ and\ s_{\beta} - y + z = t_{\gamma'}\]
then call it \textit{3-transposition}. If $\beta' - \alpha' = \gamma - \beta$ and $\gamma' - \beta' = \beta - \alpha$ then call $S$ and $T$ \textit{relative}.

\item[(T2)]
The transposition of this type operates on 2 parts. \\
For parts $\alpha, \beta \in \{1, \dots, k\}$ $(\alpha < \beta)$ on $S$, and $0 \leq x \leq s_{\alpha}$,
\[ S = (10^{s_{\alpha} - x}[0^{x}][(10^{\ast})^{\beta - \alpha - 1}10^{s_{\beta}}](10^{\ast})^{\alpha - \beta - 1 + k}) \]
is transformed into
\[ S' = (10^{s_{\alpha} - x}(10^{\ast})^{\beta - \alpha - 1}10^{s_{\beta} + x}(10^{\ast})^{\alpha - \beta - 1 + k}). \]
If there exists $\alpha', \beta' \in \{1, \dots, k\}$ on $T$ such that
\[s_{\alpha} - x = t_{\alpha'} \ and\ s_{\beta} + x = t_{\beta'}\]
then call it \textit{2-transposition}. If $\beta' - \alpha' = \beta - \alpha$ then call $S$ and $T$ \textit{relative}.
\end{enumerate}

%\label{}
%\newtheorem{thm}{Theorem}
%\newtheorem{lem}[thm]{Lemma}
%\newtheorem{cor}[thm]{Corollary}
%\newdefinition{rmk}{Remark}
%\newproof{pf}{Proof}

\noindent\textbf{Observation 1.}\label{triangle}
For $S = (0^{s_{1}}1 \dots 0^{s_{k}}1)$, $T = (0^{t_{1}}1 \dots 0^{t_{k}}1)$, $1 \leq \alpha < \beta < \gamma \leq k$ and $1 \leq \alpha' < \beta' < \gamma' \leq k$, there is 3-transposition if and only if $s_{\alpha} + s_{\beta} + s_{\gamma} = t_{\alpha'} + t_{\beta'} + t_{\gamma'}$, $t_{\alpha'} \leq s_{\alpha} + s_{\beta}$, $t_{\beta'} \leq s_{\gamma} + s_{\alpha}$ and $t_{\gamma'} \leq s_{\beta} + s_{\gamma}$. \\

\noindent Proof. 
The given conditions of 3-transposition follow the definition. Conversely, let $t_{\alpha'} > s_{\alpha} + s_{\beta}$ for any $\alpha, \beta$ and $\alpha'$, then this does not form 3-transposition. The contradiction follows for other forms in the same way. 
$\square$\\

For $S = (0^{s_{1}}1 \dots 0^{s_{k}}1)$ and $1 \leq r \leq k$, we define the following functions.
\[ f_1(S, r) = max\{s_{i_{1}} + s_{i_{2}} + \dots + s_{i_{r}} | 1 \leq i_{1} < i_{2} < \dots < i_{r} \leq k\} \]
\[ f_2(S, r) = min\{s_{i_{1}} + s_{i_{2}} + \dots + s_{i_{r}} | 1 \leq i_{1} < i_{2} < \dots < i_{r} \leq k\} \]
For $r > k$, set $f_1(S, r) = f_1(S, k)$ and $f_2(S, r) = f_2(S, k)$. Exceptionally for $r \leq 0$, set $f_1(S, r) = f_2(S, r) = 0$. \\

\noindent\textbf{Theorem 2.}\label{newfunction}
For circular binary strings $S$ and $T$, if $d(S, T) \leq m$ then
\[ \forall r,\ f_2(S, r - m) \leq f_2(T, r) \leq f_1(T, r) \leq f_1(S, r + m). \] 

\noindent Proof. 
This is shown by induction on $m$.
In the case of $m = 1$, $S$ is transformed into $T$ by relative 3-transposition or 2-transposition. Set $S = (0^{s_{1}}1 \dots 0^{s_{k}}1)$ and $T = (0^{t_{1}}1 \dots 0^{t_{k}}1)$. Let $X = \{ 1, \dots, k \}$.
\begin{enumerate}
\item[(i)]
In the case of 3-transposition, from Observation 1, for $1 \leq \alpha < \beta < \gamma \leq k$ and $1 \leq \alpha' < \beta' < \gamma' \leq k$, it forms that $s_{\alpha} + s_{\beta} + s_{\gamma} = t_{\alpha'} + t_{\beta'} + t_{\gamma'}$,  $t_{\alpha'} \leq s_{\alpha} + s_{\beta}$, $t_{\beta'} \leq s_{\gamma} + s_{\alpha}$ and $t_{\gamma'} \leq s_{\beta} + s_{\gamma}$. Let $I \subseteq X$ such that $|I| = r$ arbitrarily. Set $\sum_{i \in I}t_{i} = e(T, I)$. There are 4 cases with respect to whether $\alpha'$, $\beta'$ and $\gamma'$ are the elements of $I$ or not. In the case of $\alpha', \beta', \gamma' \notin I$, each $t_{i}$ for $i \in I$ is equal to some $s_{j}$ for $j \in X$, that leads $e(T, I) \leq f_1(S, r) \leq f_1(S, r + 1)$. In the case of $\alpha' \in I$ and $\beta', \gamma' \notin I$, $t_{\alpha'} \leq s_{\alpha} + s_{\beta}$, that leads $e(T, I) \leq f_1(S, r + 1)$. In the case of $\alpha', \beta' \in I$ and $\gamma' \notin I$, $t_{\alpha'} + t_{\beta'} \leq s_{\alpha} + s_{\beta} + s_{\gamma}$, that leads $e(T, I) \leq f_1(S, r + 1)$. In the case of $\alpha', \beta', \gamma' \in I$, $t_{\alpha'} + t_{\beta'} + t_{\gamma'} = s_{\alpha} + s_{\beta} + s_{\gamma}$, that leads $e(T, I) \leq f_1(S, r) \leq f_1(S, r + 1)$. For other combinations, the inequality holds similarly.

\item[(ii)]
In the case of 2-transposition, for $1 \leq \alpha < \beta \leq k$ and $1 \leq \alpha' < \beta' \leq k$, it forms that $s_{\alpha} + s_{\beta} = t_{\alpha'} + t_{\beta'}$, $t_{\alpha'} \leq s_{\alpha}$ and $t_{\beta'} \leq s_{\alpha} + s_{\beta}$. Let $I \subseteq X$ such that $|I| = r$ arbitrarily. Set $\sum_{i \in I}t_{i} = e(T, I)$. There are 4 cases with respect to whether $\alpha'$ and $\beta'$ are the elements of $I$ or not. In the case of $\alpha', \beta' \notin I$, each $t_{i}$ for $i \in I$ is equal to some $s_{j}$ for $j \in X$, that leads $e(T, I) \leq f_1(S, r) \leq f_1(S, r + 1)$. In the case of $\alpha' \in I$ and $\beta' \notin I$, $t_{\alpha'} \leq s_{\alpha}$, that leads $e(T, I) \leq f_1(S, r) \leq f_1(S, r + 1)$. In the case of $\beta' \in I$ and $\alpha' \notin I$, $t_{\beta'} \leq s_{\alpha} + s_{\beta}$, that leads $e(T, I) \leq f_1(S, r + 1)$. In the case of $\alpha', \beta' \in I$, $t_{\alpha'} + t_{\beta'} = s_{\alpha} + s_{\beta}$, that leads $e(T, I) \leq f_1(S, r) \leq f_1(S, r + 1)$. 
\end{enumerate}

In the case of $m \geq 2$, let $S'$ be the first transformation from $S$ by 1 transposition. $S'$ is transformed into $T$ by $m-1$ transpositions and suppose $f_1(T, r) \leq f_1(S', r + m - 1)$. $S$ is transformed into $S'$ by 1 transposition so that $f_1(S', r + m - 1) \leq f_1(S, r + m)$. We obtain $f_1(T, r) \leq f_1(S, r + m)$. For $r > k$, it forms from the definition.

$f_2(S, r - m) \leq f_2(T, r)$ is led in the same way.
$\square$\\

\noindent\textbf{Lemma 3.}\label{upper bound}
For $S = (0^{s_{1}}1 \dots 0^{s_{k}}1)$, $T = (0^{t_{1}}1 \dots 0^{t_{k}}1)$, and $X = \{1, \cdots, k\}$, let for some $l \in X$ and for all $i \in X$, $t_l \geq s_i$. If for some $i \in X$, $s_i \leq \Sigma_{j \in X - \{l\}} t_j$ then $d(S, T) \leq k - 2$. \\

\noindent Proof. 
It is obvious for $k = 2$. Let $k \geq 3$. Let $L = \{ i \in X | s_i \leq \Sigma_{j \in X - \{l\}} t_j \} \ne \emptyset$. Suppose for all $i \in L$, $s_i < \min_{j \in X} t_j$. Take $\alpha \in L, \beta \notin L$, $\gamma \ne \alpha, \beta$ and $\alpha < \beta < \gamma$. Let $t_l = t_{\beta'}, t_{\beta' + \alpha - \beta} = t_{\alpha'}$ and $t_{\beta' + \alpha - \gamma} = t_{\gamma'}$. We act transpositions on $S$ and $T$. \\

\noindent ($\diamondsuit$)

\noindent Step (i) Let $1 \leq p \leq \beta' - \alpha' - 1$ (mod $k$). If $s_{\alpha + p} > t_{\alpha' + p}$ then from the side of $S$ act 
\[(\cdots 10^{t_{\alpha' + p}}[0^{s_{\alpha + p} - t_{\alpha' + p}}][1 \cdots 1]0^{s_{\gamma}}1 \cdots).\] 
Label $s_{\alpha + p} - t_{\alpha' + p} + s_{\gamma}$ as $s_{\gamma}$. If $s_{\alpha + p} < t_{\alpha' + p}$ then from the side of $T$ act 
\[(\cdots 10^{s_{\alpha + p}}[0^{t_{\alpha' + p} - s_{\alpha + p}}][1 \cdots 1]0^{t_{\beta'}}1 \cdots).\] 
Label $t_{\alpha' + p} - s_{\alpha + p} + t_{\beta'}$ as $t_{\beta'}$.

\noindent Step (ii) Let $1 \leq p \leq \alpha' - \gamma' - 1$ (mod $k$). If $s_{\beta + p} > t_{\gamma' + p}$ then from the side of $S$ act 
\[(\cdots 10^{t_{\gamma' + p}}[0^{s_{\beta + p} - t_{\gamma' + p}}][1 \cdots 1]0^{s_{\gamma}}1 \cdots).\] 
Label $s_{\beta + p} - t_{\gamma' + p} + s_{\gamma}$ as $s_{\gamma}$. If $s_{\beta + p} < t_{\gamma' + p}$ then from the side of $T$ act 
\[(\cdots 10^{s_{\beta + p}}[0^{t_{\gamma' + p} - s_{\beta + p}}][1 \cdots 1]0^{t_{\beta'}}1 \cdots).\] 
Label $t_{\gamma' + p} - s_{\beta + p} + t_{\beta'}$ as $t_{\beta'}$.

\noindent Step (iii) Let $1 \leq p \leq \gamma' - \beta' - 1$ (mod $k$). If $s_{\gamma + p} > t_{\beta' + p}$ then from the side of $S$ act
\[(\cdots 10^{t_{\beta' + p}}[0^{s_{\gamma + p} - t_{\beta' + p}}][1 \cdots 1]0^{s_{\gamma}}1 \cdots).\]
Label $s_{\gamma + p} - t_{\beta' + p} + s_{\gamma}$ as $s_{\gamma}$. If $s_{\gamma + p} < t_{\beta' + p}$ then from the side of $T$ act
\[(\cdots 10^{s_{\gamma + p}}[0^{t_{\beta' + p} - s_{\gamma + p}}][1 \cdots 1]0^{t_{\beta'}}1 \cdots).\]
Label $t_{\beta' + p} - s_{\gamma + p} + t_{\beta'}$ as $t_{\beta'}$. \\

\noindent $S$ and $T$ was transformed into $U$ and $V$ respectively. $d(S, T) \leq d(U, V) + k - 3$. It forms $s_{\alpha} + s_{\beta} + s_{\gamma} = t_{\alpha'} + t_{\beta'} + t_{\gamma'}$, $s_{\alpha} \leq t_{\alpha'} + t_{\gamma'}$, $s_{\beta} \leq t_{\beta'} + t_{\gamma'}$ and $t_{\gamma'} \leq s_{\alpha} + s_{\beta}$. From Observation 1 by relative 3-transposition, $d(U, V) = 1$. That is, $d(S, T) \leq 1 + k - 3 = k - 2$. \\

\noindent Suppose for some $\alpha \in L$, $\min_{j \in X} t_j \leq s_{\alpha} \leq \min_{j \in X} t_j + \max_{j \in X - \{l\}} t_j$. Let $\max_{j \in X - \{l\}} t_j = t_{\alpha'}, t_l = t_{\beta'}, \min_{j \in X} t_j = t_{\gamma'}, s_{\alpha + \beta' - \alpha'} = s_{\beta}, s_{\alpha + \beta' - \gamma'} = s_{\gamma}$ and $\alpha' < \beta' < \gamma'$. We act transpositions on $S$ and $T$ as Case (*). $S$ and $T$ was transformed into $U$ and $V$ respectively. $d(S, T) \leq d(U, V) + k - 3$. It forms $s_{\alpha} + s_{\beta} + s_{\gamma} = t_{\alpha'} + t_{\beta'} + t_{\gamma'}$, $s_{\alpha} \leq t_{\alpha'} + t_{\gamma'}$, $s_{\beta} \leq t_{\beta'} + t_{\gamma'}$ and $t_{\gamma'} \leq s_{\alpha} + s_{\beta}$. From Observation 1 by relative 3-transposition, $d(U, V) = 1$. That is, $d(S, T) \leq 1 + k - 3 = k - 2$. \\

\noindent Otherwise, for some $\alpha \in L$, $s_{\alpha} > \min_{j \in X} t_j + \max_{j \in X - \{l\}} t_j$. Let $\max_{j \in X - \{l\}} t_j = t_{\alpha'}, t_l = t_{\beta'}, \min_{j \in X} t_j = t_{\gamma'}, s_{\alpha + \beta' - \alpha'} = s_{\beta}, s_{\alpha + \beta' - \gamma'} = s_{\gamma}$ and $\alpha' < \beta' < \gamma'$. We label the next as $a_1, \cdots, a_{k - 3}$.
\[ t_{\alpha' - 1} = a_1, \cdots, t_{\gamma' + 1} = a_{\alpha' - \gamma' - 1 + k}, \]
\[ t_{\beta' - 1} = a_{\alpha' - \gamma' + k}, \cdots, t_{\alpha' + 1} = a_{\beta' - \gamma' - 2 + k}, \]
\[ t_{\beta' + 1} = a_{\beta' - \gamma' - 1 + k}, \cdots, t_{\gamma' - 1} = a_{k - 3}. \]
Take minimum $m$ such that
\[ s_{\alpha} - (a_1 + \cdots + a_m) \leq t_{\alpha'} + t_{\gamma'}. \]
We act transpositions on $S$ and $T$. \\

\noindent Step 1.  Act transpositions on $S$. \\
\noindent Case (i) Let $1 \leq m \leq \alpha' - \gamma' - 1 + k$. For $0 \leq p \leq q - 2$ and $q = m$, act 
\begin{equation}
(\cdots [1][0^{*}\Sigma_{j = 0}^{p}(10^{a_{q - j}})]0^{s_{\alpha} - \Sigma_{j = q - p}^{q} a_j}1 \cdots). 
\end{equation}
For $p = m - 1$ and $r = \beta' - \gamma' + k$ act 
\begin{equation}
(\cdots [(10^{*})^{r}1][0^{*}\Sigma_{j = 0}^{p}(10^{a_{p + 1 - j}})]0^{s_{\alpha} - \Sigma_{j = 1}^{p + 1} a_j}1 \cdots). 
\end{equation}

\noindent Case (ii) Let $\alpha' - \gamma' + k \leq m \leq \beta' - \gamma' - 2 + k$. For $0 \leq p \leq \alpha' - \gamma' - 3 + k$ and $q = \alpha' - \gamma' - 1 + k$, act (1). For $p = \alpha' - \gamma' - 2 + k$ and $r = \beta' - \gamma' + k$, act (2).
For $\alpha' - \gamma' - 1 + k \leq p \leq m - 2$ and $r = \beta' - \gamma' + k$, act 
\begin{equation}
(\cdots \Sigma_{j = 0}^{\alpha' - \gamma' - 2 + k}(10^{a_{\alpha' - \gamma' - 1 + k - j}})(10^{*})^{r}10^{s_{\alpha} - \Sigma_{j = 1}^{p + 1} a_j}[\Sigma_{j = 0}^{p - (\alpha' - \gamma' - 1 + k)}(0^{a_{p + 1 - j}}1)0^{*}][1] \cdots). 
\end{equation}
Let $\Sigma_{j = 0}^{\alpha' - \gamma' - 2 + k}(10^{a_{\alpha' - \gamma' - 1 + k - j}})(10^{*})^{\beta' - \gamma' + k} = S_0$.
For $p = m - 1$, act
\[ (\cdots S_010^{s_{\alpha} - \Sigma_{j = 1}^{p + 1} a_j}[\Sigma_{j = 0}^{p - (\alpha' - \gamma' - 1 + k)}(0^{a_{p + 1 - j}}1)0^{*}][(10^{*})^{\beta' - \gamma' - p - 3 + k}1] \cdots). \]

\noindent Case (iii) Let $\beta' - \gamma' - 1 + k \leq m \leq k - 3$. For $0 \leq p \leq \alpha' - \gamma' - 3 + k$ and $q = \alpha' - \gamma' - 1 + k$, act (1). For $p = \alpha' - \gamma' - 2 + k$ and $r = 1$, act (2). For $\alpha' - \gamma' - 1 + k \leq p \leq \beta' - \gamma' - 3 + k$ and $r = 1$, act (3). Let $\Sigma_{j = 0}^{\alpha' - \gamma' - 2 + k}(10^{a_{\alpha' - \gamma' - 1 + k - j}})10^{*} = S_1$ and $\Sigma_{j = 0}^{\beta' - \alpha' - 2}(0^{a_{\beta' - \gamma' - 2 + k - j}}1) = S_2$.
For $\beta' - \gamma' - 2 + k \leq p \leq m - 2$, act
\[ (\cdots [1][0^{*}S_1\Sigma_{j = 1}^{p + 1 - (\beta' - \gamma' - 2 + k)}(10^{a_{\beta' - \gamma' - 2 + k + j}})]0^{s_{\alpha} - \Sigma_{j = 1}^{p + 1} a_j}1S_2 \cdots). \]
For $p = m - 1$, act
\[ (\cdots [(10^{*})^{k - p - 4}1][0^{*}S_1\Sigma_{j = 1}^{p + 1 - (\beta' - \gamma' - 2 + k)}(10^{a_{\beta' - \gamma' - 2 + k + j}})]0^{s_{\alpha} - \Sigma_{j = 1}^{p + 1} a_j}1S_2 \cdots). \]

\noindent Step 2. Label $s_{\alpha} - \Sigma_{j = 1}^{m} a_j$ as $s_{\alpha}$. We act transpositions on $S$ and $T$ as ($\diamondsuit$). \\

\noindent After the transpositions of Step 1 and 2, $S$ and $T$ was transformed into $U$ and $V$ respectively. $d(S, T) \leq d(U, V) + k - 3$. It forms $s_{\alpha} + s_{\beta} + s_{\gamma} = t_{\alpha'} + t_{\beta'} + t_{\gamma'}, s_{\alpha} \leq t_{\alpha'} + t_{\gamma'}$ and $s_{\beta} \leq t_{\beta'} + t_{\gamma'}$. Suppose $s_{\alpha} < t_{\gamma'}$ then $s_{\alpha} + a_m < t_{\gamma'} + a_m \leq t_{\alpha'} + t_{\gamma'}$ that contradicts the minimality of $m$. So $t_{\gamma'} \leq s_{\alpha} \leq s_{\alpha} + s_{\beta}$. From Observation 1 by relative 3-transposition, $d(U, V) = 1$. That is, $d(S, T) \leq 1 + k - 3 = k - 2$. $\square$ \\

\noindent The maximum transposition distance on strings is called \textit{transposition diameter}. Theorem 6 shows a necessary and sufficient condition of transposition diameter on circular binary strings. \\\\

\noindent\textbf{Theorem 4.}\label{diameter}
For circular binary strings $S$ and $T$ of $k \geq 2$ parts, $d(S, T) = k - 1$ if and only if $f_1(S, 1) > f_1(T, 1) \geq f_2(T, 1) > f_2(S, k - 1)$. \\

\noindent Proof. 
It is assumed that $f_1(T, 1) \leq f_1(S, 1)$. For $k = 2$ it is obvious. Let $k \geq 3$. From Theorem 2, if $f_2(T, 1) > f_2(S, k - 1)$ then $d(S, T) \geq k - 1$, and $d(S, T) \leq k - 1$ necessarily, so $d(S, T) = k - 1$. Conversely, let $f_2(T, 1) \leq f_2(S, k - 1)$, and $d(S, T) \leq k - 2$ is shown by Lemma 3. 
$\square$

\section{Conclusions}
\ \ \ \ \ This study took up transposition distances on circular binary strings. In this paper, the distance bounds were shown in terms of partitions and the diameter was represented by a necessary and sufficient condition, that is related to a NP-hard problem. It can expand to resolve transposition distances on arbitrary strings and the rest of decision problems. They could be generalized with combinatorial analysis for perspective. 

\paragraph{Acknowledgement}
This study is supported by Ota and Oda laboratory at Keio University. I appreciate their general advice of mathematics, especially for Prof. Katsuhiro Ota. 

%\label{}



%% References


%%\bibliographystyle{elsarticle-num}
%%\bibliography{<your-bib-database>}



\end{document}